\documentclass[letterpaper, 10 pt, conference]{ieeeconf}
\IEEEoverridecommandlockouts                              
\usepackage{cite}
\usepackage{amsmath,amssymb,amsfonts}
\usepackage{algorithmic}
\usepackage{graphicx}
\usepackage{textcomp}
\usepackage{amsmath,amssymb,color}
\usepackage{graphicx,epsfig,framed}
\usepackage{mathptmx,times} 
\usepackage{epstopdf,textcase} 
\usepackage{soul,multirow,pifont} 
\usepackage{color,alltt}          
\usepackage{enumerate,siunitx}    
\usepackage{caption}
\usepackage{subcaption}
\newtheorem{theorem}{Theorem}[section]

\newtheorem{proposition}[theorem]{Proposition}

\newcommand{\dd}  {{\rm d}\hbox{\hskip 0.5pt}}
\newcommand{\bbm}[1]{\left[\begin{matrix} #1 \end{matrix}\right]}

\title{\Large \bf Numerical and Lyapunov-Based Investigation of the Effect of Stenosis on Blood Transport Stability Using a Control-Theoretic PDE Model of Cardiovascular Flow}
\author{Shantanu Singh and Nikolaos Bekiaris-Liberis  
\thanks{Funded by the European Union (ERC, C-NORA, 101088147). Views and opinions expressed are however those of the authors only and do not necessarily reflect those of the European Union or the European Research Council. Neither the European Union nor the granting authority can be held responsible for them. }
\thanks{S.~Singh ({\tt\small shantanu.singh46@gmail.com}) and N.~Bekiaris-Liberis ({\tt\small nlimperis@tuc.gr}) are with the School of Electrical and Computer Eng., Technical University of Crete, Chania 731 00, Greece. }
}
\begin{document}
\maketitle
\thispagestyle{empty}
\pagestyle{empty}
\begin{abstract}
We perform various numerical tests to study the effect of (boundary) stenosis on blood flow stability, employing a detailed and accurate, second-order finite-volume scheme for numerically implementing a partial differential equation (PDE) model, using clinically realistic values for the artery's parameters and the blood inflow. The model consists of a baseline $2\times 2$ hetero-directional, nonlinear hyperbolic PDE system, in which, the stenosis' effect is described by a pressure drop at the outlet of an arterial segment considered. We then study the stability properties (observed in our numerical tests) of a reference trajectory, corresponding to a given time-varying inflow (e.g., a periodic trajectory with period equal to the time interval between two consecutive heartbeats) and stenosis severity, deriving the respective linearized system and constructing a Lyapunov functional. Due to the fact that the linearized system is time varying, with time-varying parameters depending on the reference trajectories themselves (that, in turn, depend in an implicit manner on the stenosis degree), which cannot be derived analytically, we verify the Lyapunov-based stability conditions obtained, numerically. Both the numerical tests and the Lyapunov-based stability analysis show that a reference trajectory is asymptotically stable with a decay rate that decreases as the stenosis severity deteriorates. 
\end{abstract}


\section{Introduction}\label{sec1}
Blood flow dynamics prediction and monitoring is of significant importance as it may enable accurate and timely detection of potential human health threats. In particular, one of the most common threat is related to arterial stenosis, see, e.g., \cite{Chen}, \cite{Quarteroni}, \cite{Stergiopulos_Thesis}. For this reason, there exist accurate PDE cardiovascular flow dynamic models, describing blood transport on its natural domain (that is continuous in time/space) and aiming at prediction and analysis of its dynamics, including the effect of a potential stenosis, see, for example, \cite{Bikia}, \cite{Canic}, \cite{Clark}, \cite{Formaggia1}, \cite{Quarteroni}, \cite{Reymond}. Due to computational complexity of detailed cardiovascular flow models (possibly evolving on varying, 3-D and 2-D domains), one may have to employ simpler, nevertheless accurate, 1-D PDE blood flow models in the presence of stenosis, primarily consisting of $2\times 2$ hetero-directional, nonlinear hyperbolic systems \cite{Koeppl}, \cite{Seeley}, \cite{Stergiopulos}, \cite{Young}. Motivated by these reasons, we employ here the model from \cite{Nikos_1} (which is essentially a modification of the model in \cite{Stergiopulos}, to be recast as a control-theoretic model, with a specific, modified formulation for the right boundary condition to capture the effect of stenosis), for performing numerical investigations and theoretically analyzing the effect of stenosis on blood flow stability.

Although, to the best of our knowledge, we are not aware of a control-theoretic formulation of a $1$-D PDE-based model of cardiovascular flow in the presence of stenosis apart from \cite{Nikos_1}, there exist, accurate 1-D PDE models of blood flow. Here we review the ones that are most closely related to the model we employ, which incorporate a baseline model, consisting of a second-order system of nonlinear hyperbolic PDEs (mainly describing mass/momentum conservation), properly modified to account for the different phenomena studied. These include, for example, the study of the effect of stenosis \cite{Koeppl}, \cite{Seeley}, \cite{Stergiopulos}, \cite{Young} and prosthetics \cite{Formaggia1}, as well as the incorporation in the model of heart dynamics \cite{Formaggia3} and dynamics due to other parts of the arterial network \cite{Reymond}.  To the best of our knowledge, a PDE-based, control-theoretic stability analysis of blood flow in the presence of stenosis has not been conducted in existing literature.

In the present paper, we use the model from \cite{Nikos_1}, consisting of a $1$-D hetero-directional, nonlinear hyperbolic PDE system, with the right boundary condition formulated to capture the effect of a stenosis; essentially, considering a flow bottleneck at the right boundary of an artery segment considered. We present new numerical investigations implementing  a detailed and accurate second-order, finite-volume numerical scheme (as compared to the potentially diffusive, finite difference-based numerical scheme in \cite{Nikos_1}), for performing various numerical tests to study the effect of the boundary stenosis on blood flow. In particular, we numerically solve the model for obtaining the flow and cross-sectional area, with clinically realistic parameters and inflow, of an artery segment corresponding to a part of abdominal aorta, for various degrees of stenosis severity.

 Moreover, we study the (open-loop) stability of a reference trajectory, corresponding to a given inflow and stenosis degree, based on the time-varying linearized PDE system.  Specifically, we construct a Lyapunov functional and derive the respective stability conditions. Given that the linearized system is time-varying, involving the reference trajectories themselves, which correspond to the solutions of the nonlinear hyperbolic system (for given parameters and time-varying inflow), which are not available analytically, the obtained stability conditions cannot be verified analytically. Thus, the ($L^2$) stability conditions obtained from utilization of a Lyapunov functional are verified numerically for various levels of stenosis severity. In general, a given reference trajectory, corresponding to a given inflow and parameters (including the stenosis degree) is shown to be asymptotically stable, with the stability properties, in particular, the decay rate, deteriorating with increasing stenosis severity.   
%
%

\section{1-D PDE model of blood flow with stenosis}\label{sec2}
Following \cite{Nikos_1}, the 1-D approximation of cardiovascular flow dynamics is given by the following PDE 
\begin{equation}\label{NonlinearBF_A}
\hspace{-1.6cm}\frac{\partial A(x,t)}{\partial t}=-V(x,t)\frac{\partial A(x,t)}{\partial x}-A(x,t)\frac{\partial V(x,t)}{\partial x}, 
\end{equation}
\begin{equation}\label{NonlinearBF_V}
\frac{\partial V(x,t)}{\partial t}=-V(x,t)\frac{\partial V(x,t)}{\partial x}-\frac{1}{\rho}\frac{\partial P(A(x,t))}{\partial x}-K_r\frac{V(x,t)}{A(x,t)}, 
\end{equation}
where $t\geq 0$ is time, $x\in [0, L]$ is the spatial variable, $A>0$ is the cross-sectional area of the artery, $V>0$ is the average blood speed, $\rho>0$ is the blood density, $K_r>0$ is the friction parameter related to blood viscosity and $P(A)\in\mathbb{R}$ is the pressure. The pressure function $P$ can be described by 
\begin{equation}\label{Pressure}
P(A)=\frac{\beta}{A_0}\left(\sqrt{A}-\sqrt{A_0}\right),\vspace{-1mm}
\end{equation} 
where $A_0$ is the reference arterial section area at rest and $\beta=hE\sqrt{\pi}b$, where $h>0$ is the artery wall thickness, $E>0$ is Young's modulus, and $b$ is a positive parameter. 

The boundary condition at the left end is as follows  
\begin{equation}\label{BC_BF}
A(0,t)V(0,t)=Q_{\rm in}(t), 
\end{equation}
where $Q_{\rm in}>0$ is the flow at the inlet of the artery segment considered. The boundary condition on the right end (which is the location of stenosis) is obtained from the pressure drop $\Delta P$ at $x=L$ as follows (see \cite{Nikos_1})\vspace{-2mm}
$$\hspace{-2cm}\frac{\beta}{A_0}\left(\sqrt{A(L,t)}-\sqrt{A_0}\right)-R_TA(L,t)V(L,t)\vspace{-2mm}$$
\begin{equation}\label{BC_D1}
\hspace{2cm}-V(L,t)^2\frac{K_s}{\rho}\left(\frac{A(L,t)}{A_{\rm s}}-1\right)^2=0,  
\end{equation}
where the parameter $K_s>0$ is a constant and $A_{\rm s}>0$ is the cross-sectional area of the stenosis. The pressure at the right side of the stenosis location is given by the product of the total terminal resistance $R_T$ and the flow $Q(L,t)=A(L,t)V(L,t)$, where $R_T$ is a parameter related to the conditions assumed downstream of the stenosis, see e.g., \cite{Stergiopulos}. The pressure drop $\Delta P$ (consisting of the first two terms in \eqref{BC_D1}) is obtained by taking the difference of the pressure on either side of the stenosis location, i.e., based on the instantaneous cross-sectional area $A(L,t)$ and stenosis area $A_{\rm s}$. As compared with \cite{Stergiopulos}, the pressure drop here is modelled as being dependent on the ratio $A(L,t)/A_{\rm s}$ rather than $A_0/A_{\rm s}$, considering the actual cross-sectional area $A(L,t)$ immediately  before the stenosis, rather than considering the fixed reference area $A_0$.
 
We focus our attention on the sub-critical regime (that is realistic in physiological conditions, see, for example, \cite{Quarteroni}). Therefore, we restrict the domain over which $A, V$ evolve to the nonempty, connected open subset $\Omega\subset \mathbb{R}^2$, such that  
\begin{equation}\label{Omega}
\Omega=\left\{\left.\bbm{\tilde{A}\\\tilde{V}}\in\mathbb{R}^2\right| 0<\tilde{A}, \ 0<\tilde{V},\  \tilde{V}<\frac{1}{2}\sqrt{\frac{2\beta}{\rho A_0}}\tilde{ A}^{\frac{1}{4}} \right\}, 
\end{equation}
and the eigenvalues of system  \eqref{NonlinearBF_A} and \eqref{NonlinearBF_V}, given by 
\begin{equation}\label{eigen_1}
\hspace{-0.2cm}\lambda_1(A,V)=V+\frac{1}{2}\sqrt{\frac{2\beta}{\rho A_0}} A^{\frac{1}{4}},\ \lambda_2(A,V)=V-\frac{1}{2}\sqrt{\frac{2\beta}{\rho A_0}} A^{\frac{1}{4}}, 
\end{equation}
satisfy  $\lambda_2(A,V)<0<\lambda_1(A,V)$. 

We denote by $u, v$ the following Riemann coordinates 
\begin{equation}\label{Trans_Rei_u}
u=V+2\sqrt{\frac{2\beta}{\rho A_0}}A^{\frac{1}{4}},\quad v=V-2\sqrt{\frac{2\beta}{\rho A_0}}A^{\frac{1}{4}}.
\end{equation}
The inverse transformations are 
\begin{equation}\label{InverseTrans_Rei_A}
A(u,v)=\frac{\rho^2 A^2_0}{4^5 \beta^2}(u-v)^4,\quad V(u,v)=\frac{(u+v)}{2}. 
\end{equation}
Using the Riemann coordinates $u, v$ and denoting $Y=[u\ v]^\top$, we obtain the following transformed blood flow model 
\begin{equation}\label{quasi}
\frac{\partial }{\partial t}\bbm{u(x,t)\\ v(x,t)}+\mathcal{F}(Y(x,t))\frac{\partial }{\partial x}\bbm{u(x,t)\\ v(x,t)}+\mathcal{G}(Y(x,t))=0,
\end{equation}
where 
\begin{equation}\label{FY}
\hspace{-1.4cm}\mathcal{F}(Y(x,t))=\bbm{\frac{5u(x,t)+3v(x,t)}{8} & 0\\ 0 & \frac{3u(x,t)+5v(x,t)}{8}}, 
\end{equation}
\begin{equation}\label{GY}
\mathcal{G}(Y)=\bbm{f_1(u,v) \\ f_1(u,v)},\quad  f_1(u,v)=\kappa\frac{u+v}{(u-v)^4}, \vspace{-1mm}
\end{equation}
and $\kappa=\frac{4^{\frac{9}{2}}K_r\beta^2}{\rho^2A_0^2}$. 
The boundary conditions are  \vspace{-1mm}
\begin{equation}\label{BC1_Rei}
\frac{\rho^2 A^2_0}{4^{\frac{11}{2}\beta^2}}(u(0,t)+v(0,t))(u(0,t)-v(0,t))^4=Q_{\rm in}(t), \vspace{-1mm}
\end{equation}
\begin{equation}\label{BC2_Rei}
\hspace{2cm}G\left(u(L,t),v(L,t)\right)=0, \vspace{-2mm}
\end{equation}
where 
$$
\hspace{-1.8cm}G\left(u,v\right)=\rho(u-v)^2-\frac{32 \beta}{\sqrt{A_0}}-d_1(u-v)^4(u+v) \vspace{-1mm}$$
\begin{equation}\label{G}
-4K_s\rho(u+v)^2\left(d_2(u-v)^4-1\right)^2,
\end{equation}
\begin{equation}\label{d_1}
\hspace{-2.2cm}d_1=\frac{R_T\rho^2 A_0^2}{4^3\beta^2},\quad\quad d_2=\frac{\rho^2 A_0^2}{4^5\beta^2 A_{\rm s}}. \vspace{-1mm}
\end{equation}



\section{Numerical Analysis of the blood flow model}\label{sec3}
The numerical implementation is performed using a second-order, finite-volume scheme, see, for instance, \cite{Leveque, well_balanced}. Such schemes utilise the conservative properties of the nonlinear blood flow system. Hence, we represent the cardiovascular blood flow PDEs in conservative form, or, in other words, in terms of the state variables $[A(x,t)\ Q(x,t)]^\top$, where $Q(x,t)$ is the flow. The nonlinear PDE in the state variables $[A(x,t)\ Q(x,t)]^\top$ is given as follows 
\begin{equation}\label{QandA}
\frac{\partial}{\partial t}\bbm{A(x,t)\\Q(x,t)}+\frac{\partial}{\partial x}F(A(x,t),Q(x,t))+S(A(x,t),Q(x,t))=0. 
\end{equation}
In PDE \eqref{QandA}, the source term $S$ is $S(A,Q)=\bbm{0 & K_r\frac{Q}{A}}^\top$ and the flux term $F$ is $F(A,Q)=\bbm{Q & \ \frac{Q^2}{A}+\frac{\beta}{3\rho A_0}(A^{3/2})}^\top$.

The boundary condition at $x=0$ and at $x=L$ can be obtained by substituting $V=Q/A$ in equations \eqref{BC_BF} and \eqref{BC_D1}, respectively. The spatial interval $[0 \ L]$ is divided into $n=80$ elements of equal length. The time step is $\Delta t=10^{-6}~sec$. The values of $F(A,Q)$ at the interfaces of each cell are obtained using Harten, Lax, and van Leer (HLL) flux scheme, see, for instance \cite{well_balanced}.  We use the Van-Leer slope limiter for the second-order spatial discretisation as compared to diffusive slope limiters in \cite{well_balanced}. The friction term $S(A,Q)$ is updated at each cell using the semi-implicit (SI) treatment method, see \cite[Section 3.2]{well_balanced}. 

For the simulation we consider parameters of an abdominal aorta artery, and hence, the length of the artery is considered to be $6~cm$ (see, e.g., Table 1 in \cite{Stergiopulos}), while the reference radius of the artery is considered to be $r_0=0.55~cm$. The radius of the abdominal aorta considered varies from $0.58~cm$ to $0.55~cm$ over its length, therefore, for simplicity we consider the radius of the artery at rest, constant and equal to $r_0$. The radius $r_L^+$, which is the radius of the section right after the artery, i.e., $A_{\rm s}=\pi (r_L^+)^2$, is varied from $0.55~cm$ to $0.15~cm$, based on the severity of the stenosis in the simulation scenarios considered. Note that when the radius of the section right after the boundary $x=L$  is equal to the reference radius, i.e., $r_L^+=r_0$, then this implies that there is no stenosis (in other words,  when $A_{\rm s}=A_0=\pi r_0^2$). The density of blood is $\rho=1060~Kg~ m^{-3}$, the blood viscosity  $\nu=0.0035~Pa~sec$, the friction coefficient $K_r=8\pi \nu$, and  $\beta=hE\sqrt{\pi}b$, where blood vessel thickness  is $h=0.05~cm$, Young's Modulus $E=4\times 10^5 N/m^2$, constant $b=4/3$, and $K_s=1.52$. These parameters are taken from \cite{Stergiopulos_Thesis} and \cite{Seeley}. We consider the total terminal resistance $R_T=1.33\times 10^8~N ~sec~m^{-5}$, however,  other values of $R_T$ can be considered depending on the flow conditions assumed downstream of the stenosis. We refer to, e.g., \cite[Table 2]{Stergiopulos} for the range of $R_T$. 

The numerical analysis is carried out considering the inlet flow $Q_{\rm in}(t)$ to be a periodic function which is computer generated based on the Fourier series harmonics described in \cite[Table 6.1]{Stergiopulos_Thesis}. This assumption is realistic because of the periodic nature of blood flow in humans. The time period  can be considered, for example, to be the time interval between two consecutive heartbeats (see, for example, \cite{Stergiopulos_Thesis}). The magnitude of $Q_{\rm in}(t)$ is in the range of $[1.6\times 10^{-5}, \  7\times 10^{-5}]~m^3~sec^{-1}$. 
\begin{figure}[h!] 
   \centering 
  \includegraphics[height=45mm,width=85mm]{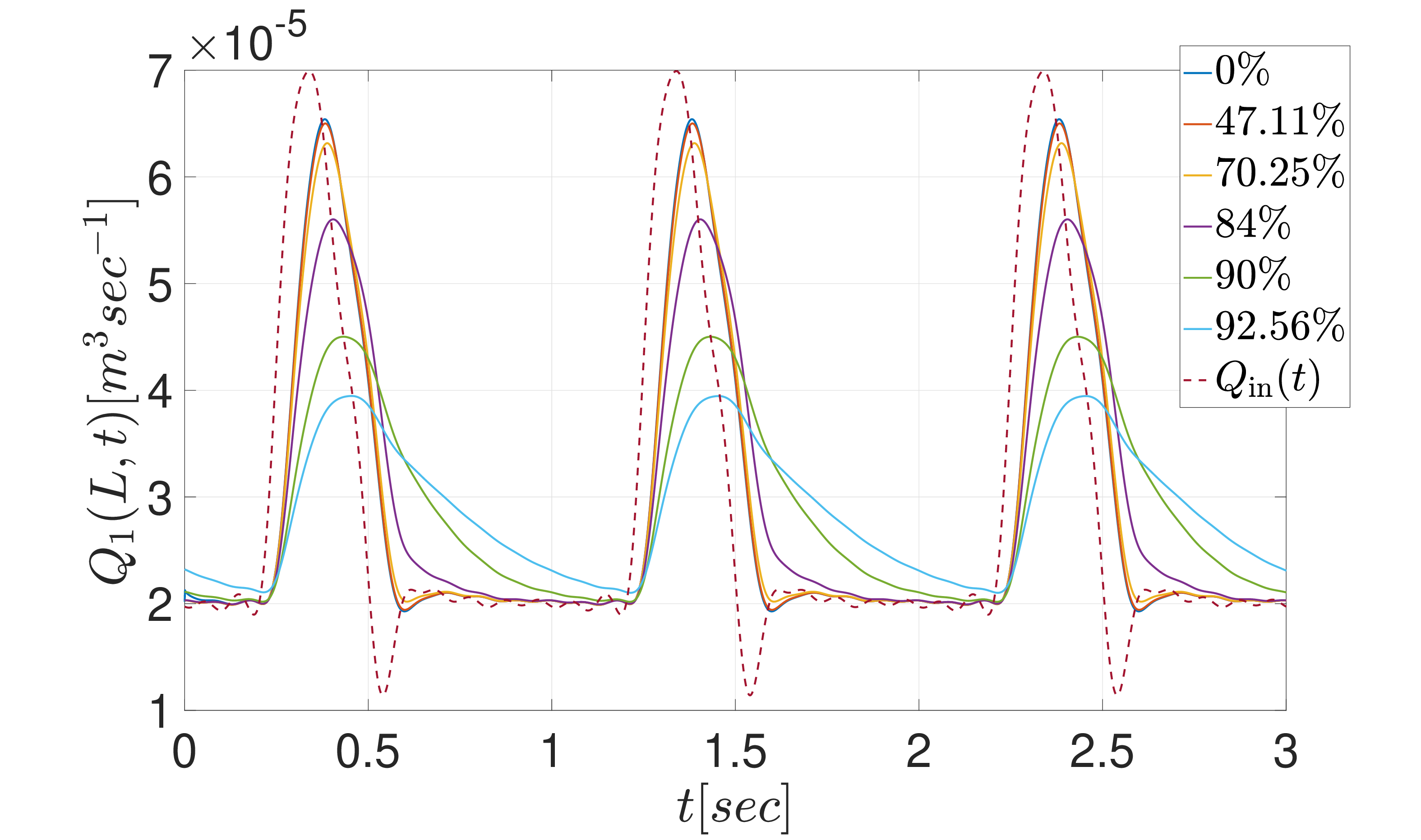}
   \caption{Flow $Q_1(L,t)$ for $0\%$, $47.11\%$, $70.25\%$, $84\%$, $90\%$,  and $92.56\%$ stenosis, computed according to $100\times\frac{A_0-A_{\rm s}}{A_0}\%$ and corresponding to $r_L^+$ values of $0.55~cm$, $0.4~cm$, $0.3~cm$, $0.22~cm$, $0.17~cm$, and $0.15~cm$, respectively.}  \label{Flow_SecondOrder_VariableStenosis}
\end{figure}
\begin{figure}[h!] 
   \centering 
  \includegraphics[height=45mm, width=85mm]{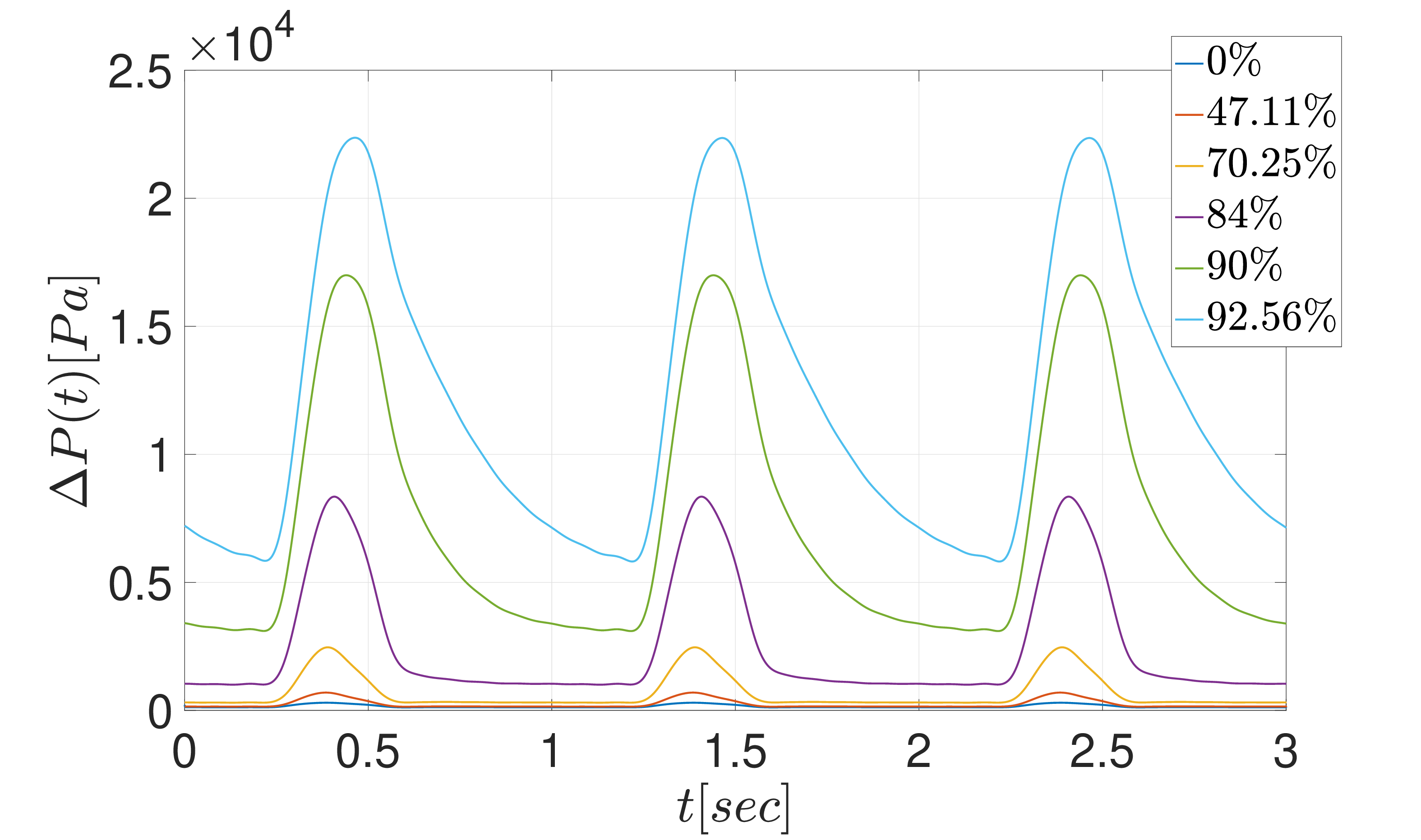}
   \caption{Pressure drop $\Delta P(t)$ for $0\%$, $47.11\%$, $70.25\%$, $84\%$, $90\%$,  and $92.56\%$ stenosis.}  \label{Pressure_SecondOrder_VariableStenosis}
\end{figure}


\begin{figure}[h!] 
   \centering 
   \includegraphics[height=45mm, width=85mm]{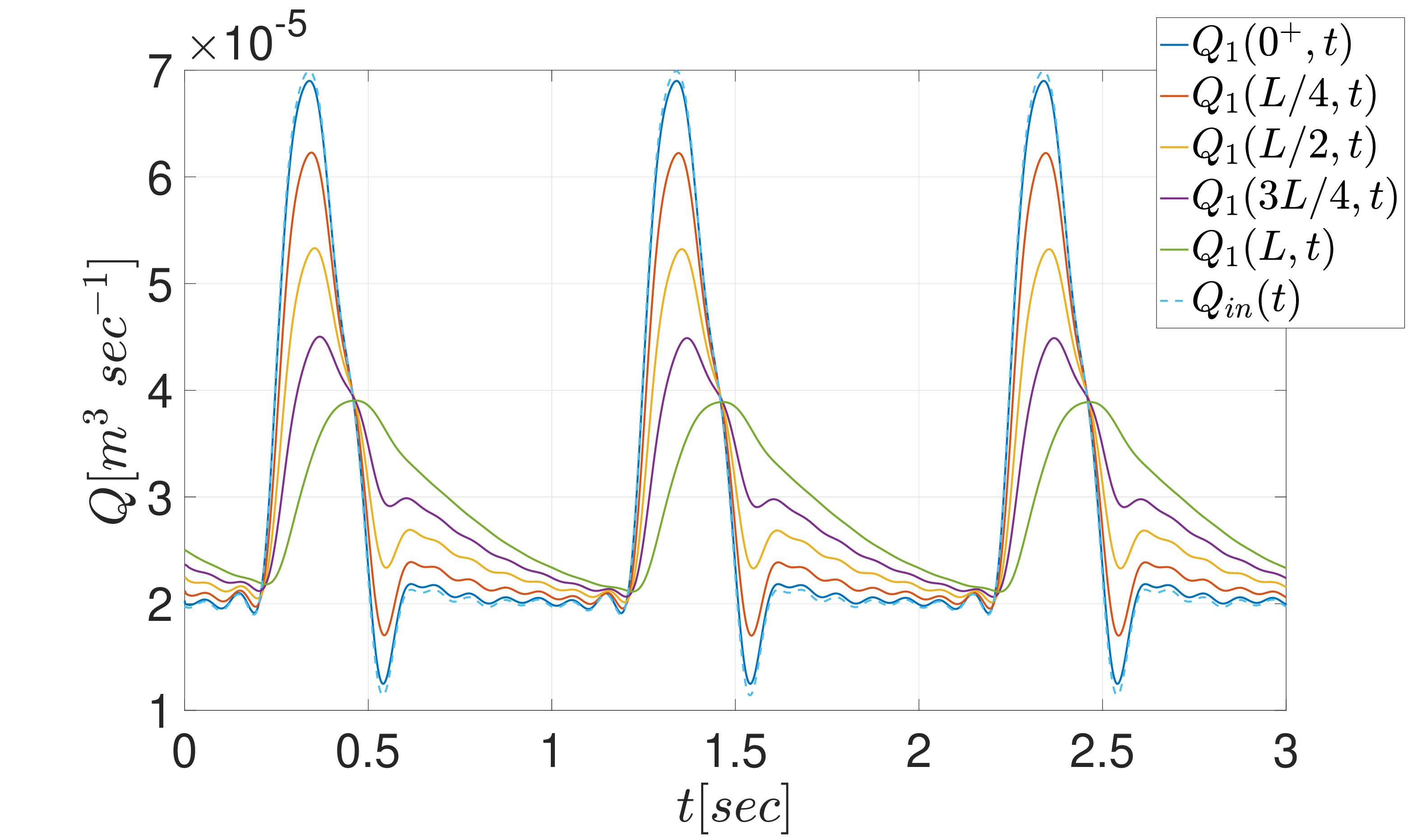}
   \caption{Flow $Q_1(x,t)$ at $x=\{\frac{L}{4}, \frac{L}{2}, \frac{3L}{4}, L\}$ when the stenosis is $92.56\%$, i.e., $r_L^+=0.15~cm$.}\label{Flow_0015}
\end{figure}
\begin{figure}[h!] 
   \centering 
   \includegraphics[height=45mm, width=85mm]{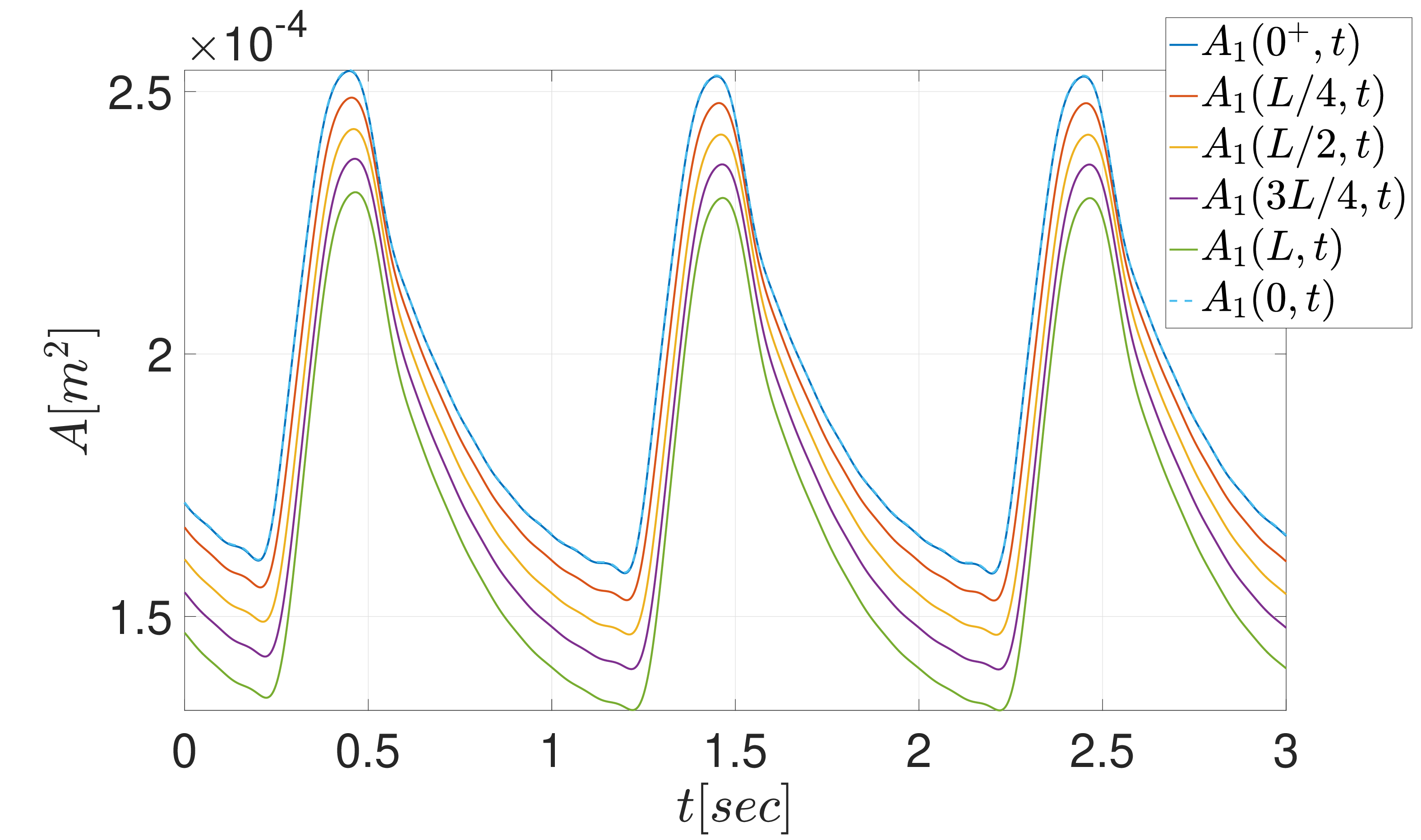}
   \caption{Area $A_1(x,t)$ at $x=\{0, \frac{L}{4}, \frac{L}{2}, \frac{3L}{4}, L\}$ when the stenosis is $92.56\%$, i.e., $r_L^+=0.15~cm$.}\label{Area_0015}
\end{figure}
%

\begin{figure}[h!] 
   \centering 
   \includegraphics[height=45mm, width=85mm]{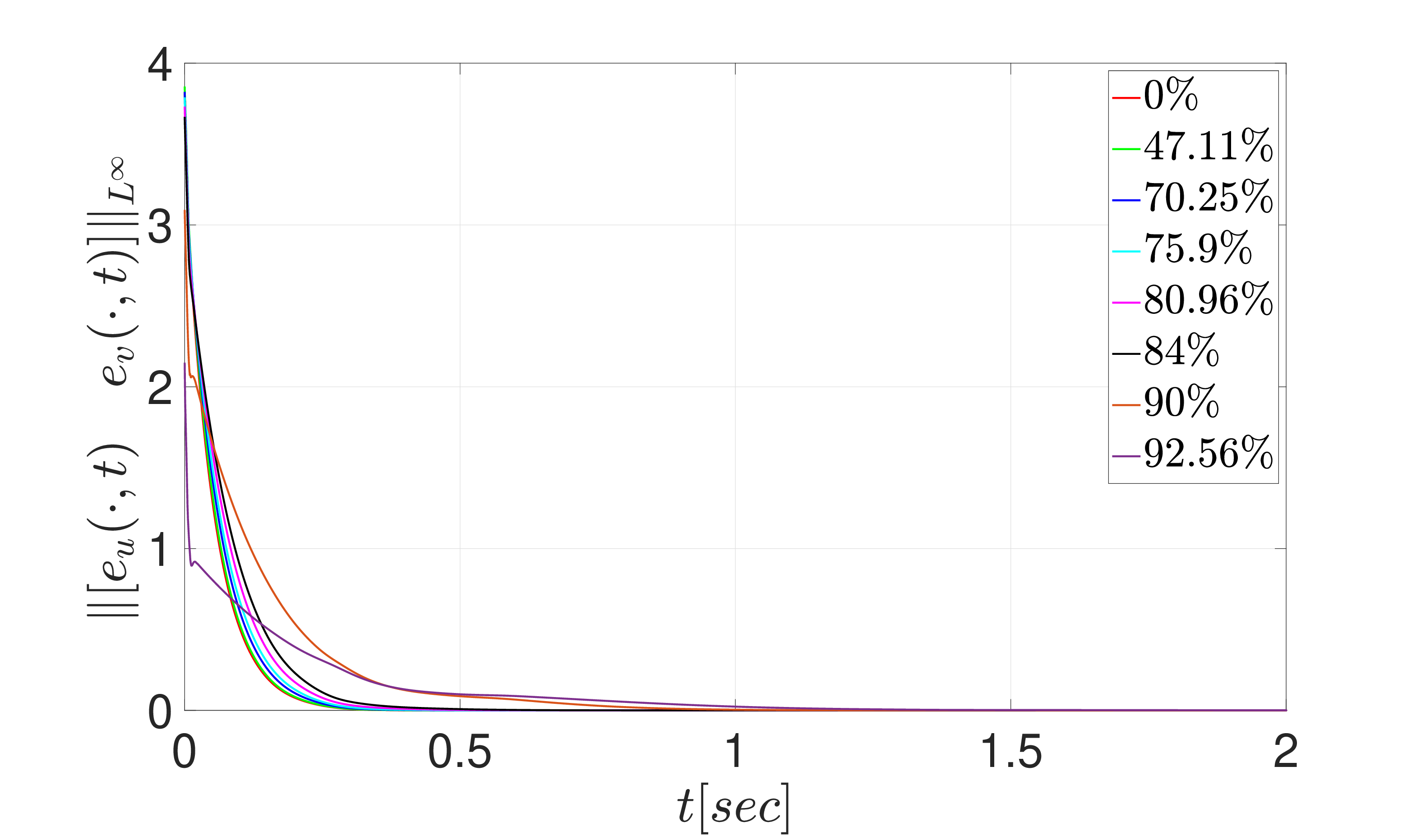}
   \caption{$L^\infty$-norm of $[e_u(\cdot,t) \ e_v(\cdot,t)]^\top$ for various levels of stenosis, where $e_u=u-u^*$ and $e_v=v-v^*$.} \label{Abs_Error_SecondOrder_VariableStenosis}
\end{figure}
The changes in flow $Q_1(L,t)$ and pressure drop $\Delta P(t)$ in the artery abdominal aorta  over time $t$  can be observed from Fig. \ref{Flow_SecondOrder_VariableStenosis} and  Fig. \ref{Pressure_SecondOrder_VariableStenosis}, respectively, with the compatible (with boundary conditions) initial conditions $[A_1(x,0) \ Q_1(x,0)]^\top$ considered to be the blood flow profile at the diastolic phase of the cardiac cycle, for each level of stenosis.  Specifically, we choose the flow/area profiles right before the start of a heartbeat, as these were observed in our simulations. As the stenosis level changes, the boundary condition at $x=L$ changes, thus, the initial condition $[A_1(x,0) \ Q_1(x,0)]^\top$ vary with stenosis severity, to be compatible with the boundary conditions. We observe that as the stenosis level increases, the peak flow $Q_1(L,t)$ during the systolic phase decreases. On the other hand, the pressure drop  $\Delta P(t)$ increases as the stenosis level increase, which is consistent with the findings, for example, in \cite{Quarteroni}, \cite{Stergiopulos}, and \cite{Stergiopulos_Thesis}. In fact, the response of the system computed via the second-order, finite-volume scheme used here, are similar to those, e.g., in \cite{Stergiopulos} and \cite{Stergiopulos_Thesis}, where finite-difference and Galerkin finite-element schemes are used.  This confirms the accuracy of the obtained numerical results.

In  Fig.~\ref{Flow_0015} and Fig. \ref{Area_0015}  we show the flow $Q_1(x,t)$ and cross-sectional area $A_1(x,t)$ at various locations, namely, $x=\{0, L/4, L/2, 3L/4, L\}$ of abdominal aorta. Notice that in Fig. \ref{Flow_0015} the magnitude of flow near $x=L$ is much lower than the flow upstream (near $x=0$). This is  due to the presence of stenosis (bottleneck) at $x=L$, the diffusive effect of the friction term, and also due to the terminal resistance ($R_T$), representing the resistance to the flow due to the effect of the downstream cardiovascular network. As $R_T$ and friction term decrease, and there is no stenosis, the flow resembles more a pure transport with delay (due to space limitation we do not include a respective plot here).

We next use initial condition $[A_2(x,0) \ Q_2(x,0)]^\top$ corresponding to the blood profile at the peak of flow during systolic phase of the cardiac cycle and in the artery with no stenosis. We keep the same initial condition $[A_2(x,0) \ Q_2(x,0)]^\top$ for all levels of stenosis to compare the error response. We observe from Fig. \ref{Abs_Error_SecondOrder_VariableStenosis} that in the case of no stenosis $[A_2\   Q_2]^\top$ trajectory converges to the reference trajectory $[A_1 \  Q_1] ^\top$ in approximately $0.5~sec$. However, when the severity of the stenosis increases to $90\%$ and $92.56\%$, then the error trajectories converge to zero in approximately $1.5~sec$. In particular, we denote by $[u^*\  v^*]^\top$ the reference Riemann coordinates corresponding to the reference state variables $[A_1\  V_1]^\top$ with the initial condition $[A_1(x,0) \ Q_1(x,0)]^\top$ (described below Fig.~\ref{Area_0015}), inflow $Q_{\rm in}(t)$, and define the error variables as $e_u=u-u^*$ and $e_v=v-v^*$, where $[u\  v]^\top$ are the Riemann coordinates that correspond to solutions originating from $[A_2(x,0) \ Q_2(x,0)]^\top$ and the same $Q_{\rm in}$.   We observe from Fig. \ref{Abs_Error_SecondOrder_VariableStenosis} that the error system is stable (this was verified also with other initial conditions, but due to space limitation we do not include here the corresponding plots), though the decay rate decreases when the stenosis degree increases. To provide a control-theoretic interpretation of these stability properties, we conduct the Lyapunov stability analysis of the error system in Section~\ref{sec5}. From a more physical viewpoint, the fact that the decay rate decreases as the stenosis degree increases can be interpreted as follows. As the stenosis becomes more severe the capacity flow at the bottleneck caused by the stenosis decreases. This in turn implies that the maximum (discharge) outflow decreases. Consequently, the transient dynamics of the respective trajectories over the domain of the arterial segment considered become slower.  


\section{Lyapunov stability of the linearized system}\label{sec5}
We analyze the nonlinear system \eqref{NonlinearBF_A} and \eqref{NonlinearBF_V} in the vicinity of a reference state trajectory, via the linearized system.  
In fact, the reference state trajectory is considered to be the one originating from the initial condition $[A_1(x,0)\  V_1(x,0)]^\top$ and corresponding to the inflow $Q_{\rm in}(t)$ (given in  Section \ref{sec3}). 
Upon linearizing \eqref{quasi} about a reference solution $Y^*(x,t)=[u^*(x,t)\ v^*(x,t)]^\top$ and substituting $u(x,t)=u^*(x,t)+e_u(x,t)$ and $v(x,t)=v^*(x,t)+e_v(x,t)$, we obtain the following linear hyperbolic PDEs in the error variables $e_u(x,t)$ and $e_v(x,t)$ \vspace{-1mm}
\begin{equation}\label{Linearized_e}
\frac{\partial}{\partial t}\bbm{e_u(x,t)\\e_v(x,t)}+\Lambda(x,t)\frac{\partial}{\partial x}\bbm{e_u(x,t)\\e_v(x,t)}+\Gamma(x,t)\bbm{e_u(x,t)\\e_v(x,t)}=0, \vspace{-1mm}
\end{equation} 
where  \vspace{-1mm}
\begin{equation}\label{A}
\hspace{-0.1cm}\Lambda(x,t)=\bbm{\frac{5u^*(x,t)+3v^*(x,t)}{8} & 0\\ 0 & \frac{3u^*(x,t)+5v^*(x,t)}{8}}, \vspace{-1mm}
\end{equation}
$$\hspace{-5cm}\Gamma(x,t)=$$
\begin{equation}\label{B}
\bbm{\frac{5}{8} u^*_x(x,t)-\kappa\frac{3u^*(x,t)+5v^*(x,t)}{(u^*(x,t)-v^*(x,t))^5} & \frac{3}{8}u^*_x (x,t)+\kappa\frac{5u^*(x,t)+3v^*(x,t)}{(u^*(x,t)-v^*(x,t))^5}\\ \frac{3}{8}v_x^*(x,t)-\kappa\frac{3u^*(x,t)+5v^*(x,t)}{(u^*(x,t)-v^*(x,t))^5} & \frac{5}{8}v_x^*(x,t)+\kappa\frac{5u^*(x,t)+3v^*(x,t)}{(u^*(x,t)-v^*(x,t))^5}}.
\end{equation}
We denote by $g(u,v)$ the left-hand side of \eqref{BC1_Rei}. Taking the linear approximation of $g(u,v)$ in the vicinity of $[u^*(0,t)\ v^*(0,t)]^\top$ we obtain \vspace{-1mm}
\begin{equation}\label{BC1_linear}
\hspace{-0.8cm}e_u(0,t)=-a(t)e_v(0,t), \vspace{-1mm}
\end{equation}
where
\begin{equation}\label{a_t}
a(t)=\left.\frac{\frac{\partial}{\partial v}g(u,v)}{\frac{\partial}{\partial u}g(u,v)}\right|_{(u^*(0,t),v^*(0,t))}=-\frac{3u^*(0,t)+5v^*(0,t)}{5u^*(0,t)+3v^*(0,t)}.
\end{equation}
(Note that for $[A \ Q]^\top\in\Omega$ the term $\frac{\partial}{\partial u}g(u,v)$ is positive.) Similarly, taking the linear approximation of $G(u,v)$ (in \eqref{BC2_Rei}) in the vicinity of $[u^*(L,t)\ v^*(L,t)]^\top$ we get \vspace{-1mm}
\begin{equation}\label{BC2_linear}
e_v(L,t)=-b(t)e_u(L,t), \vspace{-1mm}
\end{equation}
where\vspace{-1mm}
\begin{equation}\label{b_t}
\hspace{2.7cm}b(t)=\left.\frac{\frac{\partial}{\partial u}G(u,v)}{\frac{\partial}{\partial v}G(u,v)}\right|_{(u^*(L,t),v^*(L,t))}.\footnote{We verified via the numerical analysis that $\frac{\partial G}{\partial v}(u^*(L,t),v^*(L,t))\neq 0$ for the parameters of the artery and stenosis severities considered in the paper.}
\end{equation}

\subsection{Stability Analysis of Error Eystem}
We investigate stability of the time-varying, linear hyperbolic PDEs  \eqref{Linearized_e}  with boundary conditions \eqref{BC1_linear} and \eqref{BC2_linear}.
\begin{proposition}\label{LyapunovTheo}
Consider system \eqref{Linearized_e}, with  \eqref{A}, \eqref{B}, and boundary conditions given by \eqref{BC1_linear} and \eqref{BC2_linear}. Assume that  $\frac{\partial}{\partial v}G(u^*,v^*)\neq 0$ and $ [u^*\ v^*] \in C^1([0,L]\times[0,+\infty); \mathbb{R}^2)$ is such that $[A^*\  V^*] \in \Omega$ and $[u^*\ v^*]$, $[u^*_x\ v^*_x] $ are uniformly bounded\footnote{Although here this is assumed, the existence and uniqueness of global, bounded $C^1$ solutions can be studied using, for e.g., \cite{Canic} and \cite[Chap. 5]{Li}.}. The zero solution $ [e_u\ e_v] ^\top=0$ is exponentially stable in $L_2$ norm provided that there exist positive constants $p_1$, $p_2$, $\delta$, and $\mu$ such that for all $t\geq 0$  and $x\in[0, L]$ the following hold \vspace{-2mm}
\begin{equation}\label{Lya_condition1}
\hspace{1cm}\frac{|\Lambda_2(0,t)|}{\Lambda_1(0,t)}\leq \frac{p_2}{p_1}, \quad \vspace{-1mm}
\end{equation}
\begin{equation}\label{Lya_condition2}
\hspace{-2cm}e^{2\mu L}\times\frac{ p_2|\Lambda_2(L,t)|}{p_1\Lambda_1(L,t)}\times b^2(t)\leq 1, \vspace{-1mm}
\end{equation}
\begin{equation}\label{Rmatrix}
\hspace{-2.4cm}\mathcal{P}\Gamma+\Gamma^\top\mathcal{P}-\frac{\dd \mathcal{P}}{\dd x}\Lambda-\mathcal{P}\frac{\partial \Lambda}{\partial x }\geq \delta I,
\end{equation}
where  $\mathcal{P}(x)={\rm diag}\{p_1e^{-\mu x},  p_2e^{\mu x}\}$.
\end{proposition}
{\it Proof.} We consider the following Lyapunov functional candidate $\mathcal{V}(z(t))$ of the error variable $z(\cdot,t)=[e_u(\cdot, t )\ e_v(\cdot, t)]^\top$ (see, e.g., \cite{Bastin_Coron} and \cite{Hu}),
\begin{equation}\label{Lyapunov}
\mathcal{V} (z(\cdot,t))=\int^L_0z^\top (x,t)\mathcal{P}(x) z(x,t) \dd x, \vspace{-1mm}
\end{equation}
where positive definite matrix $\mathcal{P}(x)$ is as in the Proposition \ref{LyapunovTheo}.
The Lyapunov functional \eqref{Lyapunov} satisfies
\begin{equation}\label{bounded_V}
m\|z\|_{L^2}\leq \mathcal{V}(z)\leq M\|z\|_{L^2},
\end{equation}
where $m={\rm min}_{x\in[0,L]}\lambda_{\rm min}(\mathcal{P}(x))$ and $M={\rm max}_{x\in[0,L]}\lambda_{\rm max}(\mathcal{P}(x))$.  Differentiating $\mathcal{V}$ with respect to $t$ along \eqref{Linearized_e} we get\vspace{-1mm}
$$\hspace{-0.9cm}\frac{\dd \mathcal{V}}{\dd t}=2\int^L_0 z^\top(x,t)\mathcal{P}(x) \left(-\Lambda(x,t)\frac{\partial z(x,t)}{\partial x}\right)\dd x \vspace{-1mm}$$
\begin{equation}\label{dotV2}
-\int^L_0 z^\top(x,t)\left(\mathcal{P} (x)\Gamma(x,t) +\Gamma^\top(x,t)\mathcal{P}(x)\right) z(x,t)  \dd x.
\end{equation}
Integrating by parts the first term in \eqref{dotV2} we obtain
$$\hspace{-2cm}\frac{\dd \mathcal{V}}{\dd t}=-[z^\top (x,t)\mathcal{P}(x) \Lambda(x,t)z(x,t)]^L_0\ +\vspace{-1mm}$$
\begin{equation}\label{dotV3}
\int^L_0 z^\top(x,t)\left(\frac{\dd \mathcal{P}}{\dd x}\Lambda+\mathcal{P}\frac{\partial \Lambda}{\partial x}-\mathcal{P}\Gamma-\Gamma^\top\mathcal{P}\right) z(x,t) \dd x.\vspace{-1mm}
\end{equation}
Denote by $R(x,t)$ the left-hand side of inequality \eqref{Rmatrix}, i.e., \vspace{-1mm}
\begin{equation}\label{R(x,t)}
R=\mathcal{P}\Gamma+\Gamma^\top\mathcal{P}-\frac{\dd \mathcal{P}}{\dd x}\Lambda-\mathcal{P}\frac{\partial \Lambda}{\partial x }.\vspace{-1mm}
\end{equation}
 Suppose that there exists a constant $\delta>0$ such that for all $x\in [0, L]$ and $t\geq 0$ inequality \eqref{Rmatrix} holds then 
\begin{equation}\label{dotV4}
\frac{\dd \mathcal{V}}{\dd t}\leq-[z^\top (x,t)\mathcal{P}(x) \Lambda(x,t)z(x,t)]^L_0-\delta \|[e_u(\cdot, t)\ e_v(\cdot, t)]\|_{L_2}^2.\vspace{-1mm}
\end{equation}
The first term in \eqref{dotV4} can be written as follows
$$\hspace{-4cm}-[z^\top (x,t)\mathcal{P}(x) \Lambda(x,t)z(x,t)]^L_0=$$
$$\hspace{2cm}\left(\mathcal{P}_1(0)\Lambda_1(0,t)a^2(t)+\mathcal{P}_2(0)\Lambda_2(0,t)\right)e^2_v(0,t) $$
\begin{equation}\label{dotV6}
- \left(\mathcal{P}_1(L)\Lambda_1(L,t)+\Lambda_2(L,t)\mathcal{P}_2(L)b^2(t)\right)e^2_u(L,t).
\end{equation}
Substituting \eqref{a_t}  and using the fact that $a(t)=-\frac{\Lambda_2(0,t)}{\Lambda_1(0,t)}$, if \eqref{Lya_condition1} holds then we obtain that\vspace{-1mm}
 \begin{equation}\label{first_part}
 \left(p_1\Lambda_1(0,t)a^2(t)+p_2\Lambda_2(0,t)\right)e^2_v(0,t)\leq 0. 
 \end{equation}
If \eqref{Lya_condition2} holds then since $\Lambda_1>0$ and $\Lambda_2<0$, we obtain
  \begin{equation}\label{second_part}
  \left(p_1\Lambda_1(L,t)e^{-\mu L}+p_2\Lambda_2(L,t)e^{\mu L}b^2(t)\right)e^2_u(L,t)\geq 0.
   \end{equation}
Thus, using \eqref{first_part} and  \eqref{second_part} we obtain from \eqref{dotV4} \vspace{-1mm}
\begin{equation}\label{dotV8}
\frac{\dd \mathcal{V}(z(\cdot,t))}{\dd t}\leq -\delta \|(z(\cdot,t))\|^2 .  \vspace{-1mm}
\end{equation} 
Hence from \eqref{bounded_V}, system \eqref{Linearized_e} is exponentially stable. $\square$

\subsection{Verification of Stability Conditions}
\begin{figure}[h!]
\centering 
\includegraphics[height=45mm, width=85mm]{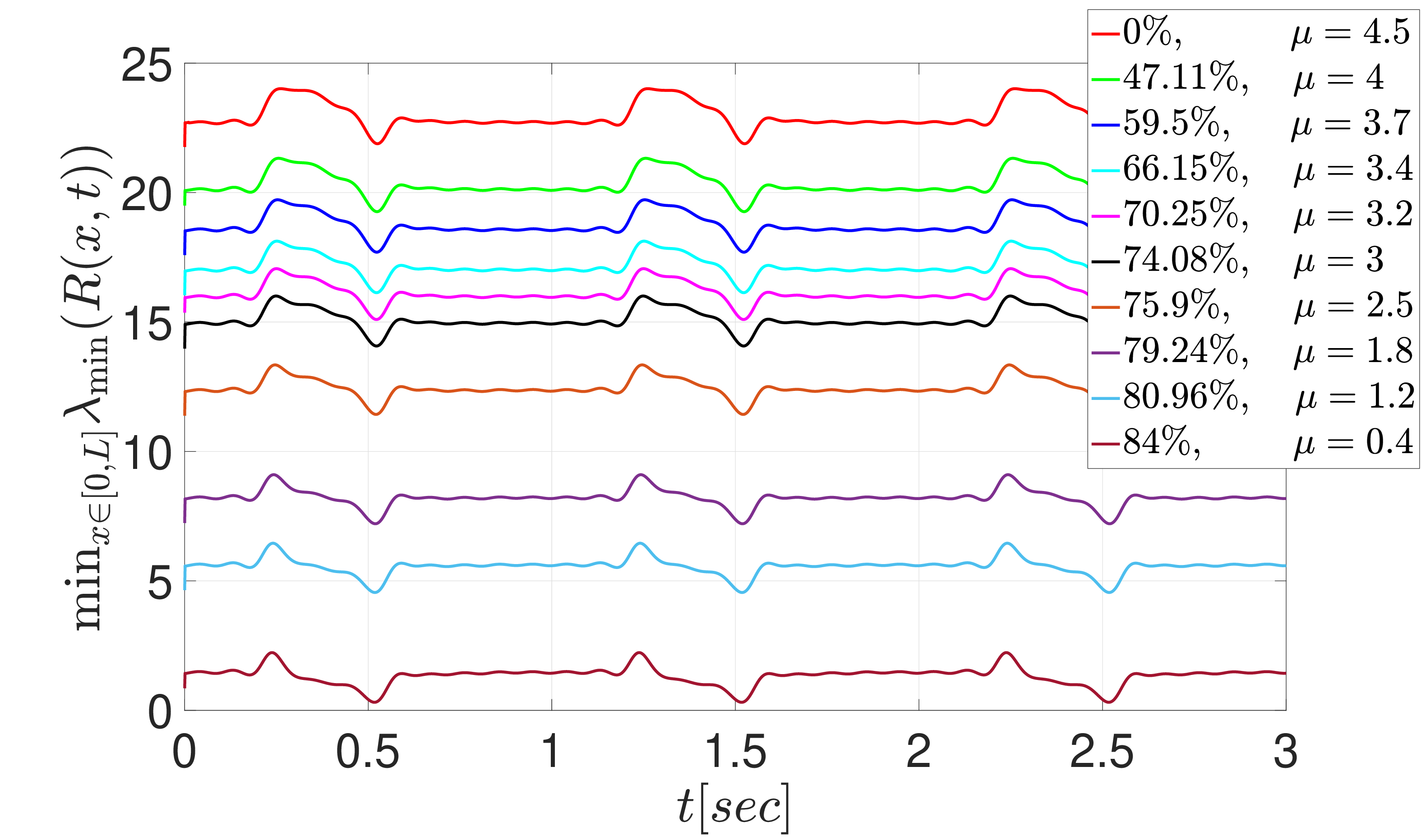}
\caption{Value of ${\rm min}_{x\in[0, L]} \lambda_{\rm min}\left(R(x,t)\right)$ for various levels of stenosis with $p_1=1$ and $p_2=0.98$.} \label{Eig_R}\vspace{-1.5mm}
\end{figure}

\begin{figure}[h!]
\centering 
\includegraphics[height=45mm, width=85mm]{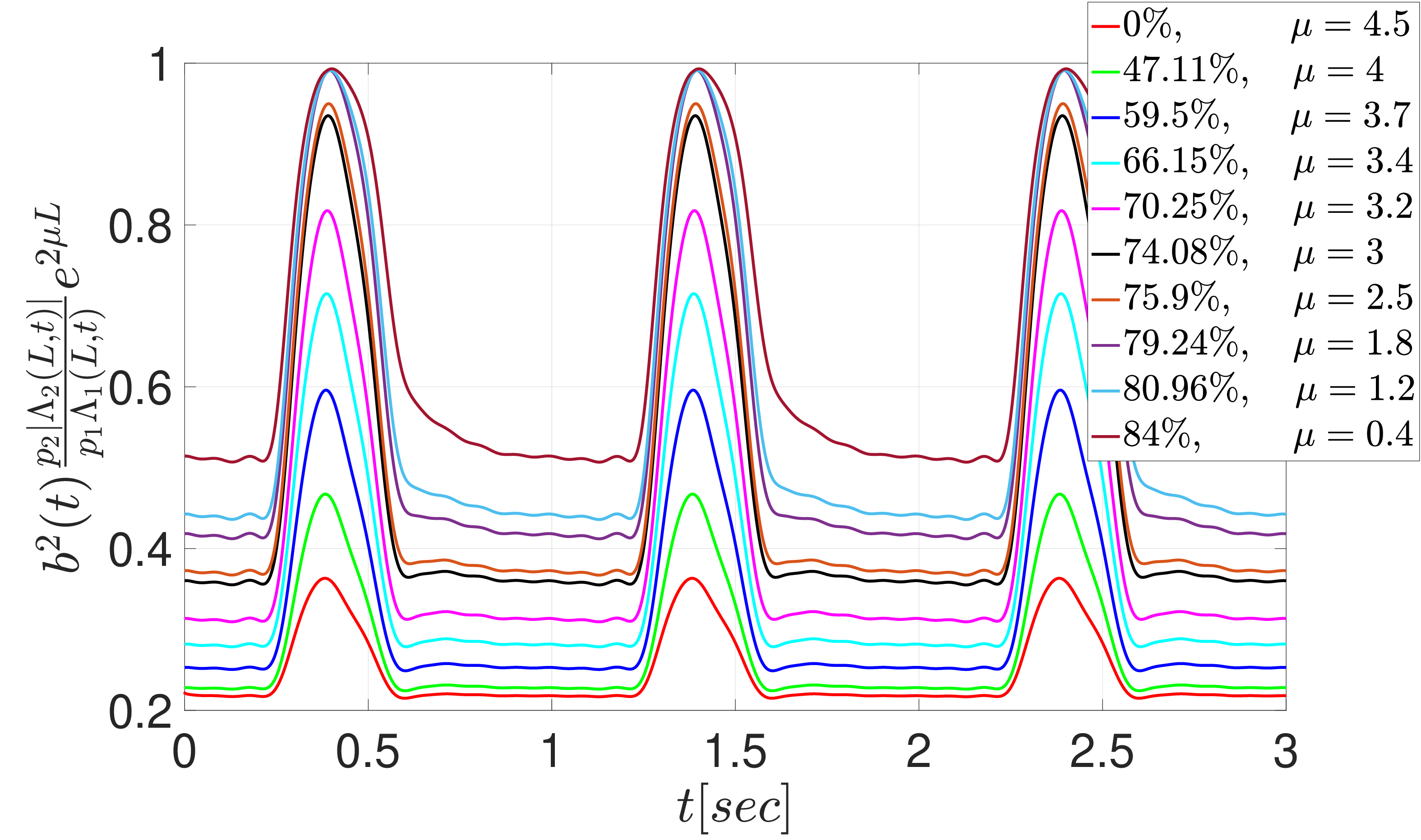}
\caption{Verification of condition \eqref{Lya_condition2} for various levels of stenosis with $p_1=1$ and $p_2=0.98$.}\label{B_cond}\vspace{-1.5mm}
\end{figure}

\begin{figure}[h!]
\centering 
\includegraphics[height=45mm, width=85mm]{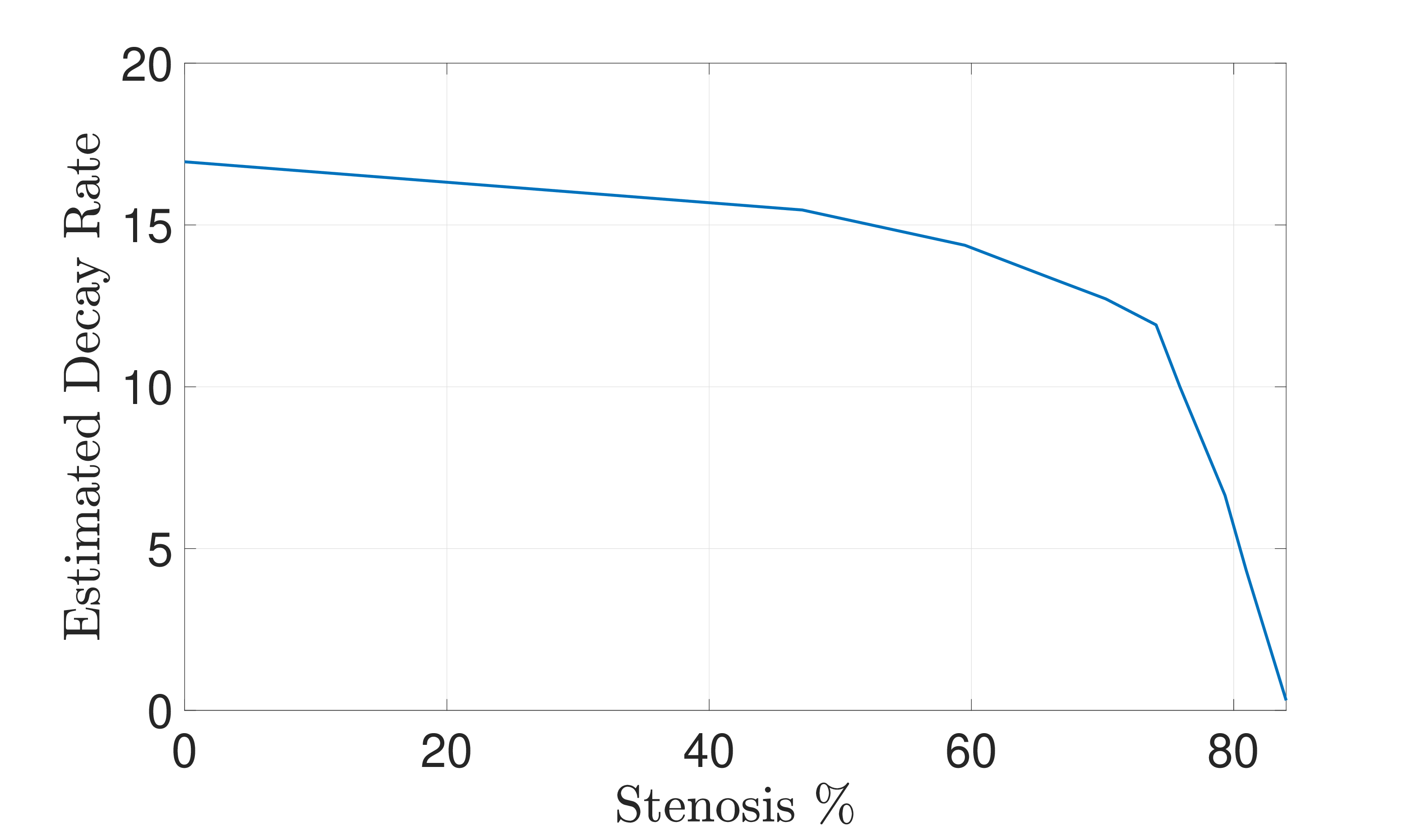}
\caption{Estimate of the decay rate given in \eqref{lambda_min} plotted against stenosis percentage.}\label{DecayRate}\vspace{-1.5mm}
\end{figure}

Now we verify the stability conditions for class of systems defined by equation \eqref{Linearized_e}. We say a ``class of systems" because for each level of stenosis severity (or each $A_{\rm s}$), the linearized system is different due to different reference trajectories $[u^*\ v^*]$ and different boundary conditions at $x=L$. For all the different cases of stenosis we select $p_1=1$ and $p_2=0.98$, in order  to satisfy inequality \eqref{Lya_condition1}\footnote{Due to space limitation we do not present a respective plot.} for all levels of stenosis severity and to also keep the ratio $p_2/p_1$ sufficiently small in order to satisfy \eqref{Lya_condition2}. We observe from Fig. \ref{Eig_R} that the minimum eigenvalues of the  matrix $R(x,t)$ are positive up to $84\%$ stenosis severity, implying that $R(x,t)$ is positive definite.  Moreover, we observe from Fig. \ref{B_cond} that inequality \eqref{Lya_condition2} is satisfied for stenosis severity of up to  $84\%$. It was observed from the numerical analysis that  the magnitude of $\mu$ should be decreased with increasing severity of stenosis, in order to satisfy \eqref{Lya_condition2}.  We observe that if for some $\delta>0$\vspace{-1mm}
\begin{equation}\label{lambda_min}
\frac{\underset{t\geq 0, x\in [0,L]}{\rm min} \lambda_{\rm min}\left(R(x,t)\right)}{M}\geq \frac{\delta}{M},\vspace{-1mm}
\end{equation}
where $M= {\rm sup}_{x\in [0, L]} \lambda_{\rm max}(\mathcal{P}(x))$, then $R(x,t)\geq \delta \mathcal{P}(x)$ for all $t\geq 0$ and $x\in [0, L]$.
Thus, the linearized error systems defined by  \eqref{Linearized_e} satisfy inequalities \eqref{Lya_condition1}--\eqref{Rmatrix} for stenosis severity of up to $84\%$. It follows from inequalities \eqref{dotV8} and \eqref{lambda_min} that system \eqref{Linearized_e} is exponentially stable with the decay rate $\delta$, having an upper bound given by the left-hand  side of \eqref{lambda_min}. The estimate of the decay rate as a function of the stenosis severity is shown in Fig. \ref{DecayRate}, from which we observe that it decreases for higher stenosis levels. This observation is consistent with the error response of the nonlinear system \eqref{QandA}, shown in Fig. \ref{Abs_Error_SecondOrder_VariableStenosis}, although in Fig.~\ref{DecayRate} there is an abrupt drop in the estimated decay rate for stenosis percentages larger than about $75\%$. This difference may be attributed to the conservatism of the Lyapunov stability conditions derived and the specific analysis parameters chosen. 

As the stenosis severity increases, the choice of positive constants $\mu$, $p_1$, and $p_2$ which simultaneously satisfy all the three inequalities \eqref{Lya_condition1}--\eqref{Rmatrix},   for all times $t\geq 0$ could not be obtained. 
For higher stenosis severity  we observed that during time periods where \eqref{Lya_condition2} was violated, matrix $R(x,t)$ had positive eigenvalues, and vice versa. This, in combination with the fact that based on the numerical simulations in Section \ref{sec3} the system is stable even for higher stenosis levels (see Fig. \ref{Abs_Error_SecondOrder_VariableStenosis}) indicate that a time-varying (possibly periodic) Lyapunov functional may be better suited for the stability analysis of our model for higher levels of stenosis. In principle, one could utilize the results from \cite{Prieur_1} to construct such a time-varying Lyapunov functional. 

\section{Conclusions}\label{sec6}
In this article we investigate the effect of stenosis on blood flow stability in abdominal aorta artery, from a control-theoretic perspective. Specifically, we investigate numerically the convergence of the trajectories of the nonlinear blood flow system to a reference trajectory. We show that the rate of convergence decreases as the stenosis severity increases. We approximate the nonlinear blood flow system by the linear system obtained via linearisation around the reference trajectory. Using Lyapunov stability analysis, we state and verify the conditions that should be satisfied for the linear, time-varying error system to be exponentially stable. In our future work, we intend to study the properties of the nonlinear system, thus also establishing the properties of the reference trajectories, which are now assumed; as well as to address the observer design problem (see, e.g., \cite{Vazquez}).

\end{document}